\documentclass[11pt]{amsart}
\usepackage{amsmath,amsfonts,amssymb}
\usepackage[all]{xy}

\makeatletter
\providecommand*{\shuffle}{%
  \mathbin{\mathpalette\shuffle@{}}%
}
\newcommand*{\shuffle@}[2]{%
  \sbox0{$#1\vcenter{}$}%
  \kern .15\ht0 
  \rlap{\vrule height .25\ht0 depth 0pt width 2.5\ht0}%
  \raise.1\ht0\hbox to 2.5\ht0{%
    \vrule height 1.75\ht0 depth -.1\ht0 width .17\ht0 %
    \hfill
    \vrule height 1.75\ht0 depth -.1\ht0 width .17\ht0 %
    \hfill
    \vrule height 1.75\ht0 depth -.1\ht0 width .17\ht0 %
  }%
  \kern .15\ht0 
}
\makeatother

\begin{document}

\newcommand{\cA}{{\mathcal A}}
\newcommand{\add}{{\rm add}}
\newcommand{\ab}{{\rm ab}}
\newcommand{\talpha}{{\tilde{\alpha}}}
\newcommand{\bB}{{\mathbb B}}
\newcommand{\cB}{{\mathcal B}}
\newcommand{\bI}{{I_\bullet}}
\newcommand{\C}{{\mathbb C}}
\newcommand{\Cat}{{\rm Cat}}
\newcommand{\CP}{{\mathbb{C}P}}
\newcommand{\cH}{{\mathcal H}}
\newcommand{\Comod}{{\rm Comod}}
\newcommand{\D}{{\mathbb D}}
\newcommand{\Diff}{{\rm Diff}}
\newcommand{\F}{{\rm F}}
\newcommand{\Fun}{{\rm Fun}}
\newcommand{\ex}{{\rm ex}}
\newcommand{\Ext}{{\rm Ext}}
\newcommand{\tE}{{\widetilde{E}}}
\newcommand{\f}{{\mathfrak f}}
\newcommand{\G}{{\mathbb G}}
\newcommand{\GT}{{\sf GT}}
\newcommand{\sG}{{\sf G}}
\newcommand{\tgamma}{{\tilde{\gamma}}}
\newcommand{\geo}{{\rm geo}}
\newcommand{\half}{{\textstyle{\frac{1}{2}}}}
\newcommand{\Hom}{{\rm Hom}}
\newcommand{\HoM}{{Ho(\mathcal M)}}
\newcommand{\loc}{{\rm loc}}
\newcommand{\bm}{{\bf m}}
\newcommand{\mm}{{\mathfrak m}}
\newcommand{\Maps}{{\rm Maps}}
\newcommand{\Mor}{{\rm Mor}}
\newcommand{\Mod}{{\rm Mod}}
\newcommand{\MTM}{{\sf MT}}
\newcommand{\cM}{{\mathcal M}}
\newcommand{\NSymm}{{\sf NSymm}}
\newcommand{\Oh}{{\rm O}}
\newcommand{\perf}{{\rm perf}}
\newcommand{\bP}{{\mathbb P}}
\newcommand{\PBW}{{\rm PBW}}
\newcommand{\Q}{{\mathbb Q}}
\newcommand{\QSymm}{{\sf QSymm}}
\newcommand{\tQSymm}{{{\bf QSymm}}}
\newcommand{\stuffle}{{\; \overline{\shuffle}\;}}
\newcommand{\R}{{\mathbb R}}
\newcommand{\Rep}{{\rm Rep}}
\newcommand{\SO}{{\rm SO}}
\newcommand{\Sp}{{\rm Sp}}
\newcommand{\Spec}{{\rm Spec}}
\newcommand{\SU}{{\rm SU}}
\newcommand{\Symm}{{\sf Symm}}
\newcommand{\TC}{{\rm TC}}
\newcommand{\THH}{{\rm THH}}
\newcommand{\tK}{{\widetilde{K}}}
\newcommand{\hK}{{\widehat{K}}}
\newcommand{\sK}{{\sf{K}_{\mathbb Q}\mathbb Z}}
\newcommand{\Tor}{{\rm Tor}}
\newcommand{\T}{{\mathbb T}}
\newcommand{\U}{{\rm U}}
\newcommand{\Wh}{{\rm Wh}}
\newcommand{\w}{{\mathfrak w}}
\newcommand{\tr}{{\rm tr}}
\newcommand{\tw}{{\widetilde{\mathfrak w}}}
\newcommand{\Z}{{\mathbb Z}}
\newcommand{\hZ}{{\widehat{\mathbb Z}}}
\newcommand{\cell}{{\rm cell}}
\newcommand{\colim}{{\rm colim}}
\newcommand{\cO}{{\mathcal O}}
\newcommand{\sS}{{s\mathbb{S}}}
\newcommand{\sOmega}{{s\Omega}}
\newcommand{\USO}{{\rm U/SO}}
\newcommand{\FO}{{\rm F/O}}
\newcommand{\fO}{{\rm f/O}}
\newcommand{\SUSO}{{\rm SU/SO}}
\newcommand{\cF}{{\mathcal F}}
\newcommand{\SpU}{{\rm Sp/U}}
\newcommand{\SpSU}{{\rm Sp/SU}}
\newcommand{\HH}{{\rm HH}}

\title{ Homotopy-theoretically enriched categories of noncommutative motives}
\author{Jack Morava}
\address{The Johns Hopkins University,
Baltimore, Maryland 21218}
\email{jack@math.jhu.edu}

\date{18 March 2015}
\begin{abstract}{Waldhausen's $K$-theory of the sphere spectrum (closely
related to the algebraic $K$-theory of the integers) is naturally
augmented as an $S^0$-algebra, and so has a Koszul dual. Classic work of Deligne
and Goncharov implies an identification of the rationalization of this 
(covariant) dual with the Hopf algebra of functions on the motivic group 
for their category of mixed Tate motives over $\Z$. This paper argues that 
the rationalizations of categories of non-commutative motives defined
recently by Blumberg, Gepner, and Tabuada consequently have natural 
enrichments, with morphism objects in the derived category of mixed Tate 
motives over $\Z$. We suggest that homotopic descent theory lifts this 
structure to define a category of motives defined not over $\Z$ but
over the sphere ring-spectrum $S^0$.} \end{abstract}

\maketitle

\section{Introduction} \bigskip

\noindent
{\bf 1.1} Building on earlier work going back at least three decades 
[26] Deligne and Goncharov have defined a $\Q$-linear Abelian
rigid tensor category of mixed Tate motives over the integers of a 
number field: in particular, the category $\MTM_\Q(\Z)$ of such motives
over the rational integers. Its generators are tensor powers $\Q(n) = 
\Q(1)^{\otimes n}$ of a Tate object $\Q(1$), inverse to the Lefschetz 
hyperplane motive (which can be regarded as a degree 
two shift of the complex
\[
0 \to \bP_1 \to \bP_0 \to 0
\]
in Voevodsky's derived category). We argue here that these
objects are analogous to the (even-dimensional) cells of stable
homotopy theory: in that, for example, the image 
\[
\bP_n = \Q(0) \oplus \cdots \oplus \Q(-n)
\]
of projective space in this category splits as a sum of terms resembling 
Lefschetz's hyperplane sections.\bigskip

\noindent
Deligne and Goncharov's definition [27 \S 1.6] depends on the validity of 
the Beilinson-Soul\'e vanishing conjecture for number fields, which implies
that their category $\MTM_\Q(\Z)$ can be characterized by a very simple 
spectral sequence with $E_2$-term
\[
\Ext^*_\MTM(\Q(0),\Q(n)) \Rightarrow K(\Z)_{2n-*} \otimes \Q
\]
equal to zero if $* > 1$ or when $* = 0, \; n \neq 0$. Borel's 
theory of regulators [19, 29] identifies the nonvanishing groups
\[
K_{4k+1}(\Z) \otimes \Q \subset \R
\]
with the subgroup of rational multiples of the conjecturally transcendental 
values $\zeta(1+2k)$ of the Riemann zeta function at odd {\bf positive} integers. 
\bigskip

\noindent
To a homotopy theorist, this is strikingly reminiscent of Atiyah's interpretation 
of Adams' work on Whitehead's homomorphism
\[
J_{n-1} : KO_n(*) = \pi_{n-1}\Oh \to \lim_{m \to \infty} \pi_{m+n-1}(S^m) = 
\pi^S_{n-1}(*) \;,
\]
ie the effect on homotopy groups of the monoid map
\[
\Oh = \lim_{m \to \infty} \Oh(m)  \to \lim_{m \to \infty} \Omega^{m-1} S^{m-1} 
:= Q(S^0) \;.
\]
The image of a Bott generator 
\[
KO_{4k}(*) \cong \Z \to {\rm image} \; J_{4k-1} \cong (\half \zeta(1-2k)\cdot \Z)/\Z \subset 
\Q/\Z
\]
under this homomorphism can be identified with the (rational) value of the zeta 
function at an odd {\bf negative} integer. \bigskip

\noindent
[Here is a quick sketch of this argument: A stable real vector bundle over 
$S^{4k}$ is classified by its equatorial twist 
\[
\xymatrix{
\talpha : S^{4k-1} \ar[r]^\alpha & \Oh \ar[r] & Q(S^0)} \;,
\]
which defines a stable cofibration
\[
\xymatrix{
S^{4k-1} \ar[r]^\talpha & S^0 \ar[r] & {\rm Cof} \; \talpha \ar[r] & 
S^{4k} \ar@{.>}[r] & \dots \;.}
\]
Adams' $e$-invariant is the class of the resulting sequence
\[
\xymatrix{
0 \ar[r] & KO(S^{4k}) \ar[r] & KO({\rm Cof \; \talpha}) \ar[r] & KO(S^0) 
\ar[r] & 0 }
\]
in a group 
\[
\Ext_{\rm Adams}^1(KO(S^0),KO(S^{4k}))
\]
of extensions of modules over the stable $KO$-cohomology operations. [$KO$ is 
contravariant, while motives are covariant, making $KO$ of a sphere analogous
to a Tate object.] After profinite completion [4] these cohomology 
operations become an action of the group $\hZ^\times$ of units in the 
profinite integers - in fact the action, through its abelianization, of the 
absolute Galois group of $\Q$ - and the resulting group of extensions can be 
calculated in terms of generalized Galois cohomology as
\[
\Ext^1_{\hZ^\times}(\hZ(0),\hZ(2k)) \cong H^1_c(\hZ^\times,\hZ(2k))
\]
(where $u \in \hZ^\times$ acts on $\hZ(n)$ as multiplication by $u^n$). At 
an odd prime $p$, $H^1_c(\hZ^\times_p,\hZ_p(2k))$ is zero unless $2k = 
(p-1)k_0$, when the group is cyclic of $p$-order $\nu_p(k_0)+1$. By 
congruences of von Staudt and Clausen, this is the $p$-ordert of the 
Bernoulli quotient
\[
\frac{B_{2k}}{2k} \in \Q/\Z \;;
\]
a global argument over $\Q$ (ie using the Chern character [1 \S 7.1b])
refines this to a homomorphism
\[
H^1_c(\hZ^\times,\hZ(2k)) \to \Q/\Z
\]
which sends a generator of $KO(S^{4k})$ to the class of $\half \zeta(1-2k)$.
See also [26 \S 3.5].]
\bigskip

\noindent
{\bf 1.2} This paper proposes an analog of the theory of mixed Tate motives
in the world of stable homotopy theory, based on B\"okstedt's theorem [[16], or
more recently [14]] that the morphism
\[
K(S^0) \to K(\Z) 
\]
of ring-spectra (induced by the Hurewicz morphism
\[
[1 : S^0 \to H\Z] \in H^0(S^0,\Z) \;)
\]
becomes an isomorphism after tensoring with $\Q$. At this point, odd zeta-values 
enter differential topology [43]. To be more precise, we argue that (unlike $K(\Z)$), 
$K(S^0)$ is naturally {\bf augmented} as a ring-spectrum over $S^0$, via the Dennis trace
\[
\tr_D : K(S^0) \to \THH(S^0) \sim S^0 
\]
[70]. Current work [38, 64] on descent in homotopy theory suggests the category of 
comodule spectra over the {\bf covariant} Koszul dual
\[
S^0 \wedge^L_{K(S^0)} S^0 := K(S^0)^\dagger 
\]
of $K(S^0)$ (or perhaps more conventionally, the category of module 
spectra over 
\[
R\Hom_{K(S^0)}(S^0,S^0) \;)
\]
as a natural candidate for a homotopy-theoretic analog of $\MTM_\Q(\Z)$. 
This paper attempts to make this plausible {\bf after tensoring $K(S^0)$ with 
the rational field} $\Q$. 
\bigskip 

\noindent
{\bf 1.3 Organization} Koszul duality is a central concern of this paper; in its most
classical form, it relates (graded) exterior and symmetric Hopf algebras. The first 
section below observes that the Hopf algebra of {\bf quasi}symmetric functions is 
similarly related to a certain odd-degree {\bf square-zero} augmented algebra. Stating
this precisely (ie over $\Z$) requires comparison of the classical shuffle product 
[31 Ch II] with the less familiar quasi-shuffle product [40, 50]. I am especially 
indebted to Andrew Baker and Birgit Richter for explaining this to me. \bigskip
 
\noindent
The next section defines topologically motivated generators (quite different from 
those of Borel) for $K_*(S^0) \otimes \Q$. Work of Hatcher [36 \S 6.4], Waldhausen [69], 
and B\"okstedt on pseudo-isotopy theory has been refined by Rognes [63] to construct an 
infinite-loop map
\[
\omega : B(\FO) \to \Wh(*) \; (\subset K(S^0))
\]
(F being the monoid of homotopy self-equivalences of the stable sphere) which
is a rational equivalence. This leads to the definition of a homotopy equivalence
\[
\w : (S^0 \vee \Sigma kO) \otimes \Q \to K(S^0) \otimes \Q
\]
of ring-spectra, with a square-zero extension of the rational sphere spectrum on 
the left, which can then be compared with Borel's calculations. Some of the work of 
Deligne and Goncharov is then summarized to construct a lift of this rational isomorphism 
to an equivalence between the algebra $\cH_{GT^*_\MTM}$ of functions on the motivic group 
of the Tannakian category $\MTM_\Q(\Z)$ and the covariant Koszul dual 
$K(S^0)^\dagger_* \otimes \Q$. \bigskip

\noindent
The final section is devoted to applications: in particular to the 
`decategorification' [18, 49 \S 4] of two-categories of `big' (noncommutative) 
motives constructed by Blumberg, Gepner, and Tabuada [10, 11], and to work of 
Kitchloo [46] on  categories of symplectic analogs of motives. The objects in 
the categories of `big' motives are themselves small stable 
$\infty$-categories, with stable  $\infty$-categories of suitably exact 
functors between them as morphism objects.  The (Waldhausen) $K$-theory 
spectra of these morphism categories define new categories enriched 
over the homotopy category of $K(S^0)$-module spectra [11 Corollary 4.13], 
having the original small stable categories as objects.\bigskip

\noindent
`Rationalizing' (tensoring the morphism objects in these homotopy categories with $\Q$) 
defines categories enriched over $K_*(S^0) \otimes \Q$-modules, to which the Koszul duality 
machinery developed here can be applied. Under suitable finiteness hypotheses, this constructs 
categories of noncommutative motives enriched over the derived category $D_b(\MTM_\Q(\Z)$ of 
classical mixed Tate motives. 
\bigskip

\noindent
{\bf Acknowledgements:} I am deeply indebted to Andrew Blumberg, Kathryn Hess,
and Nitu Kitchloo for help and ecouragement in the early stages of this work;
and to Andrew Baker, Birgit Richter, and John Rognes for their advice and
intervention in its later stages. Thanks to all of them -- and to some 
very perceptive and helpful referees -- for their interest and patience. 
The mistakes, misunderstandings, and oversimplifications below are
my responsibility.\bigskip

\section{Quasisymmetric functions and Koszul duality} \bigskip 

\noindent
{\bf 2.1} The fudamental example underlying this paper could well
have appeared in Tate's 1957 work [68] on the homology of local rings;
but as far as I know it is not in the literature, so I will begin with
it: Let $E_* := E_*[e_{2k+1} \:|\; k \geq 0]$ be the primitively generated 
graded-commutative Hopf algebra over $\Z$ with one generator in each odd 
degree, and let
\[
\phi_E : E_* \to E_*/E_+ = \Z
\]
be the quotient by its ideal $E_+$ of positive-degree elements; then 
\[
\Tor^E_*(\Z,\Z) \; \cong \; P[x_{2(k+1)} \:|\: k \geq 0] \; (:= \Symm_*)
\]
is a graded-commutative Hopf algebra with one generator in each even degree, 
canonically isomorphic to the classical algebra of symmetric functions with 
coproduct
\[
\Delta x(t) = x(t) \otimes x(t)
\]
($x(t) = \sum_{k \geq 0} x_{2k} t^k, \; x_0 := 1$). \bigskip

\noindent
This is an instance of a very general principle: if $A_* \to k$ is an augmented 
commutative graded algebra (assuming for simplicity that $k$ is a field), then 
\[
k \otimes^L_{A*} k \;= \Tor_*^A(k,k) := A_*^\dagger
\]
is an augmented, graded-commutative Hopf algebra, with 
\[
\Ext^*_A(k,k) := R\Hom^*_A(k,k) 
\]
as its graded dual [22 XVI \S 6]. More generally,
\[
(A-\Mod) \ni M \mapsto \Tor_*^A(M,k) := M^\dagger_* 
\]
extends this construction to a functor taking values in a category of graded $A^\dagger_*$-comodules.  
This fascinated John Moore [41, 60], and its implications have become quite important in 
representation theory [7, 8]; more recently, the whole subject has been vastly 
generalized by the work of Lurie.\bigskip 

\noindent
{\bf 2.2.1} For our purposes it is the quotient
\[
\varphi_\tE : E_* \to E_*/(E_+)^2 := \tE_*
\]
of the exterior algebra above, by the ideal generated by products of positive-degree elements, 
which is relevant. This quotient is the square-zero extension
\[
\tE_* = \Z \oplus \tE_+ = \Z \oplus \{e_{2k+1} \:|\: k \geq 0\}
\]
of $\Z$ by a graded module with one generator in each odd degree. \bigskip

\noindent
{\bf Proposition:} After tensoring with $\Q$, the induced homomorphism
\[
\varphi^\tE_* : \Tor^E_*(\Z,\Z) \cong \Symm_* \to \Tor^\tE_*(\Z,\Z) \cong \tQSymm_*
\]
of Hopf algebras is the inclusion of the graded algebra of rational symmetric functions 
into the algebra of rational quasi-symmetric functions, given the classical shuffle 
product $\shuffle$. \bigskip

\noindent
{\bf Proof:} In this case the classical bar resolution
\[
\bB_*(\tE/\Z) = \tE_* \otimes_\Z(\oplus_{n \geq 0} \tE_+[1]^{\otimes n})
\otimes_\Z \Z
\]
(of $\Z$ as an $\tE_*$-module [55 Ch X \S 2.3, 31 Ch II, 24 \S 2]) is, apart from 
the left-hand term, just the tensor algebra of the graded module $\tE_+[1]$ (obtained from
$\tE_+$ by shifting the degrees of its generators up by one), with algebra structure defined 
by the shuffle product; but I will defer discussing that till \S 2.3 below. Since $\tE_*$ is a 
DGA with trivial differential and trivial product, the homology $\Tor^\tE_*(\Z,\Z)$ of the complex
\[
\Z \otimes_\tE \bB_*(\tE/\Z) = \oplus_{n \geq 0}(\tE_+[1])^{\otimes n}
\]
with its resulting trivial differential is the algebra $\tQSymm_*$ on
$\tE_+[1]$. [Tate, by the way, worked with a commutative noetherian local ring    
\[
\phi_A : A \to B = A/\bm_A = k 
\]
and studied $\Tor^A_*(k,k)$, though not as a Hopf algebra; but in his calculations he 
used what is visibly the resolution above.] \bigskip      

\noindent
{\bf Remark:} In fact under very general conditions [38 \S 3.1, 39 \S 6.12] the bar
construction associates to a morphism $\varphi : A \to B$ of suitable monoid objects, a pullback 
functor
\[
L\varphi^* : M_* \mapsto M_* \otimes_A \bB_*(A/B) \cong M_* \otimes^L_A B 
\]
from some (simplicial or derived) category of modules over $\Spec \; A$ to a similar category 
of modules over $\Spec \; B$, cf \S 2.4 below. Here $\bB_*(A/B)$ is a resolution of $B$ as
an $A$-module, corresponding to $B(A,B,A)$ in [32 Prop 7.5], cf also [31, 53, \dots]. In 
the example above, regarding $E_*$ as a DGA with trivial differential, we obtain a covariant 
functor from the bounded derived category of $E_*$-modules to the bounded derived category of 
graded modules with a {\bf co}action of the classical Hopf algebra of symmetric functions. 
However, the algebra of symmetric functions is canonically self-dual over $\Z$ [54 I \S 4]) 
and we can interpret this derived pullback as a functor to the bounded derived category of 
modules over the dual symmetric algebra. \bigskip

\noindent
{\bf 2.2.2} In this paper I will follow K Hess [38 \S 2.2.23 - 2.2.28]: a morphism
\[
\varphi : A \to B
\]
of monoids in a suitable (eg simplicially enriched [38 \S 3.16, \S 5.3]) 
category of modules (perhaps over a differential graded algebra or a 
ring-spectrum) defines an $A$-bimodule {\bf bi}algebra
\[
B \wedge^L_A B \; := \; W(\varphi)
\]
(analogous to an algebraic topologist's Hopf algebroid, though in general
without antipode). In her framework the construction above is a {\bf descent} functor: 
its target has a natural (`Tannakian') enrichment [37 \S 5.3, 43] or lift $L\varphi^\dagger$ 
to a category of $B$-modules with compatible coaction by the descent coring $W(\varphi)$.
Lurie's work (eg [52 \S 7.13], [53 \S 5.2.2 - 5.2.3]) provides the natural context for such
constructions. \bigskip

\noindent
{\bf 2.2.3} Completing the argument requires clarifying relations between the shuffle product
$\shuffle$ and the quasi-shuffle or `stuffle' product $\stuffle$. A {\bf shuffle} of a pair $r,s 
\geq 1$  of integers is a partition of the set $\{1,\dots,r+s\}$ into disjoint subsets $a_1 < 
\cdots < a_r$ and $b_1 < \cdots < b_s$; such a shuffle defines a permutation 
\[
\sigma(1,\dots,r+s) = (a_1,\cdots,a_r,b_1,\dots,b_s) \;.
\]
The shuffle product on the tensor algebra $T^\bullet(V)$ of a module $V$ is 
defined by 
\[
v_1 \cdots v_r \shuffle v_{r+1} \cdots v_{r+s} = \sum v_{\sigma(1)} \cdots 
v_{\sigma(r+s)} \;,
\]
with the sum taken over all shuffles of $(r,s)$. The deconcatenation coproduct
\[
\Delta : T^\bullet(V) \to T^\bullet(V) \otimes T^\bullet(V)
\]
sends $v_1 \cdots v_r$ to the sum
\[
v_1 \cdots v_r \otimes 1 \; + \; \sum_{1< i < r} v_1 \cdots v_i \otimes 
v_{i+1} \cdots v_r \; + \;  1 \otimes v_1 \cdots v_r \;.
\]
The algebra $T^\bullet(V)$, with this (commutative but not cocommutative) Hopf structure, 
is sometimes called the cotensor (Hopf) algebra of $V$. The shuffle product is characterized 
by the identity
\[
(v \cdot x) \shuffle (w \cdot y) = v \cdot (x \shuffle (w \cdot y)) + w \cdot ((v 
\cdot x) \shuffle y) \;,
\] 
where $v,w \in V$ and $x,y \in T^\bullet(V)$. I will write $\tQSymm_*$ for the Hopf algebra
$\Tor^\tE_*(\Z,\Z)$ of \S 2.2.1, with $\shuffle$ as product. \bigskip

\noindent
The closely related Hopf algebra $\QSymm_*$ of quasi-symmetric functions over $\Z$, with the 
{\bf quasi}-shuffle product $\stuffle$ , is perhaps most efficiently defined as dual to the free graded 
associative Hopf algebra 
\[
\NSymm^* := \Z \langle Z_{2(k+1)} \:|\: k \geq 0 \rangle
\]
of noncommutative symmetric functions [6, 22 \S 4.1.F, 36], with coproduct
\[
\Delta Z(t) = Z(t) \otimes Z(t) 
\]
($Z(t) = \sum_{k \geq 0} Z_{2k} t^k, \; Z_0 = 1$).\bigskip

\noindent
More generally, if $(V,\star)$ is a (graded) commutative algebra, the 
quasi-shuffle [or overlapping shuffle, or stuffle] product $\stuffle$ on 
$T^\bullet(V)$ is a deformation [40 \S 6] of the shuffle product characterized 
by the identity
\[
(v \cdot x) \stuffle (w \cdot y) = v \cdot (x \stuffle (w \cdot y)) + w \cdot((v \cdot x) 
\stuffle y) + (v \star w) \cdot (x \stuffle y) \;.
\]
In particular, if we define an algebra structure on the graded vector space 
spanned by classes $f_i$ dual to the $Z_i$'s by $f_i \star f_j = f_{i+j}$, 
we recover the quasi-shuffle product on the dual of $\NSymm^*$. \bigskip

\noindent
The Lie algebra of primitives in $\NSymm^*$ is generated by the analogs of
Newton's power functions [23 \S 4.1.F], and map under abelianization
to the classical power function primitives in $\Symm^*$; dualizing yields
a morphism of $\Symm_*$ to $\QSymm_*$, which rationalizes to the asserted 
inclusion. $\Box$ \bigskip

\noindent
{\bf 2.3 Remarks} \bigskip

\noindent
i) The applications below will be based on a variant of $\tE_*$ defined by
generators in degree $4k+1, \; k \geq 0$, rather than $2k+1$. The corresponding
free Lie algebras will then have generators in (homological) degree $-2(2k+1)$.
This doubling of topological degree relative to motivic weight is a familiar
consequence of differing conventions. \bigskip

\noindent
ii) Hoffman [40 Theorem 2.5] constructs an isomorphism
\[
\exp : \tQSymm_* \otimes \Q \to \QSymm_* \otimes \Q
\]
of graded Hopf algebras over the rationals, taking $\shuffle$ to $\stuffle$; so 
over $\Q$ we can think of the morphism defined by the proposition as the inclusion 
of the symmetric functions in the quasiymmetric functions with the quasishuffle product. \bigskip

\noindent
iii) The rationalization $\NSymm^* \otimes \Q$ is the (primitively generated) 
universal enveloping algebra $U(\f^*)$ of the free Lie algebra $\f^*$ generated
by the $Z$'s over $\Q$. By Poincar\'e-Birkhoff-Witt its modules can 
be regarded as representations of a pro-unipotent groupscheme $\sG_0(\f^*)$ over 
$\Q$, or equivalently as comodules over the Hopf algebra $\QSymm_* \otimes \Q$ 
of algebraic functions on that pro-unipotent group. If we interpret graded modules 
as representations of the multiplicative groupscheme in the usual way [[3 \S 3.2.7], 
see also [8 \S 1.1.2]] then we can regard these modules as representations of a proalgebraic 
groupscheme
\[
1 \to \sG_0(\f^*) \to \sG(\f^*) := \G_m \ltimes \sG_0(\f^*) \to \G_m \to 1 \;. 
\]
In a very helpful appendix, Deligne and Goncharov [27 \S A.15] characterize 
representations of $\sG(\f^*)$ as graded $\f^*$-modules, such that (if 
$\Q(n)$ denotes a copy of $\Q$ in degree $n$) 
\[
\Ext^1_{\rm{Rep}(\sG(\f^*))}(\Q(0),\Q(n)) = (\f^n_\ab)^\vee \;.
\] 
This is explained in more detail in [35 \S 8]; we will return to 
this description below. \bigskip

\noindent
iv) The rational stable homotopy category is equivalent to the derived 
category of rational vector spaces, and the homotopy category of rational 
ring-spectra is equivalent to the homotopy category of DGAs: the Hurewicz
map
\[
[S^*,X_\Q] \cong \pi^S_*(X) \otimes \Q \to H_*(X,\Q)
\]
is an isomorphism. This leads to a convenient abuse of notation which may not
distinguish the rationalization $X_\Q$ of a spectrum from its homology (or its
homotopy). For example, the rational de Rham algebra of forms on a reasonable
space is a good model for the rational Spanier-Whitehead commutative 
ringspectrum $[X,S^0_\Q]$. \bigskip

\noindent
{\bf 2.4} Finite-dimensional graded modules over a field $k$ have a good 
duality functor 
\[
V_* \mapsto \Hom_k(V_*,k) = (V^*)^\vee \;,
\]
and a great deal of work on the homological algebra of augmented algebras 
$\phi : A \to k$ (and their generalizations) is formulated in terms of 
constructions generalizing
\[
M_* \to \Hom_k(M_* \otimes_A \bB_*(A/k),k) \cong \Hom_A(M_*,\bB^*(A/k)^\vee)
:= R\Hom_A(M_*,k) \;,
\]
where $\bB^*(A/k)^\vee$ is now essentially a {\bf co}bar construction. This is 
the classical contravariant Koszul duality functor: see [12, 21, 30 \S 4.22] for 
recent work in the context of modules over ring spectra. \bigskip

\noindent
In this paper, however, we work instead [following [34, 60] with the 
{\bf co}variant functor $M_* \mapsto M^\dagger_*$ defined as in \S 2.2.2, but 
regarded as mapping modules over an augmented algebra $A$ to comodules 
over a coaugmented coalgebra $A^\dagger$. In particular, in the case of our
main example (over $\Q$) Hess's hypotheses [38 \S 5.3] are statisfied, and we 
have the \bigskip

\noindent
{\bf Corollary:} The Hess-Koszul-Moore functor 
\[
L\varphi_\tE^\dagger : D_b(\tE_* \otimes \Q - \Mod) \to D_b(\tQSymm_* \otimes 
\Q - \Comod)
\]
is an equivalence of symmetric monoidal categories. \bigskip

\noindent
{\bf Proof:} The point is that, in the context of graded commutative augmented 
algebras $A$ over a field $k$, the functor $L\varphi^*$ is monoidal, in the sense
that
\[
L\varphi^*(M_0) \otimes_k L\varphi^*(M_1)  := (M_0 \otimes_A \bB(A/k)) \otimes_k 
k \otimes_k (M_1 \otimes_A \bB(A/k)) \otimes_k k 
\]
is homotopy-equivalent to 
\[
L\varphi(M_0 \otimes_A M_1) := ((M_0 \otimes_A M_1) \otimes_A \bB(A/k)) \otimes_k k 
\]
via the morphism
\[
(M_0 \otimes_A \bB(A/k)) \otimes_k M_1 \to M_0 \otimes_A M_1
\]
induced by the homotopy equivalence of $\bB(A/k)$ with $k$ as an $A$-module. This
then lifts to an equivalence $L\varphi^\dagger$ of comodules, cf [58, 68]. $\Box$ 
\bigskip

\noindent
In the terminology of \S 2.2.3iii, the composition
\[
\PBW \circ \exp^* \circ L\varphi_\tE^\dagger := L\Phi^\dagger_\tE
\]
thus defines an equivalence of the derived category of $\tE_* \otimes 
\Q$-modules with the derived category of $\sG(\f^*)$-representations. A 
similar argument identifies the bounded derived category of modules over 
$E_* \otimes \Q$ with the bounded derived category of representations of 
the graded abelianization $\sG(\f^*_\ab)$ of $\sG(\f^*)$: in other words, 
of graded modules over $\Symm^* \otimes \Q$. \bigskip

\noindent
{\bf Remarks} The quite elementary results above were inspired by 
groundbreaking work of Baker and Richter [6], who showed that the integral 
cohomology of $\Omega \Sigma \C P^\infty_+$ is isomorphic to $\QSymm^*$ as 
a Hopf algebra. Indeed, the $E_2$-term 
\[
\Tor^{H^*(Z)}(H^*(X),H^*(Y)) \Rightarrow H^*(X \times_Z Y)
\]
of the Eilenberg-Moore spectral sequence for the fiber product
\[
\xymatrix{
\Omega \Sigma \CP^\infty_+ \ar@{.>}[d] \ar@{.>}[r] & P \Sigma \CP^\infty_+ 
\ar[d] \\
{*} \ar[r] & \Sigma \CP^\infty_+ }
\]
is the homology of the bar construction on the algebra $H^*(\Sigma \C 
P^\infty_+)$, which is a square-zero extension of $\Z$. The spectral sequence 
collapses for dimensional reasons, but has nontrivial multiplicative 
extensions (connected to the fact that $H^*(\C P^\infty)$ is polynomial
(cf \S 2.2.2, [5, 50]). \bigskip

\noindent
In view of Proposition 3.2.1 below, Proposition 2.2.1 can be rephrased as the 
algebraically similar assertion that the K\"unneth spectral sequence [32 IV \S 4.1] 
\[
\Tor^{H_*(K(S^0),\Q)}(H_*(S^0,\Q),H_*(S^0,\Q)) \Rightarrow 
H_*(K(S^0)^\dagger,\Q)
\]
for ring-spectra collapses. In this case the algebra structure on 
$H_*(K(S^0),\Q)$ is trivial, resulting in the shuffle algebra $\tQSymm_*$.
Note that although these spectral sequences look algebraically similar,
one is concentrated in positive, the other in negative, degrees. If or how
they might be related, eg via the cyclotomic trace (cf \S 4.1), seems 
quite mysterious to the author.\bigskip

\noindent
I am deeply indebted to Baker and Richter for help with this, and
with many other matters. I am similarly indebted to John Rognes for
patient attempts to educate me about the issues in the
section following. \bigskip

\section{Geometric generators for $K_*(S^0)\otimes \Q$} \bigskip

\noindent
{\bf 3.1} Stable smooth cell bundles are classified by a space (ie simplicial 
set)
\[
\colim_{n \to \infty} \; B\Diff(\D^n) \;,
\]
where $\Diff(\D^n)$ is the group of diffeomorphisms of the {\bf closed}
$n$-disk (which are {\bf not} required to fix the 
boundary\begin{footnote}{This paper was inspired by Graeme Segal's 
description of such objects as `blancmanges'}\end{footnote}). Following 
[72 \S 1.2, \S 6.1], there is a fibration
\[
\cH_\Diff (S^{n-1}) \to B\Diff(\D^n) \to B\Diff(\D^{0n}) \;,
\]
where $\D^{0n}$ is the {\bf open} disk, and $\cH_\Diff(S^{n-1})$ is the 
simplicial set of smooth $h$-cobordisms of a sphere with itself [68]. The 
homomorphism $\Oh(n) \to \Diff(\D^{0n})$ is a homotopy equivalence, while 
the constructions of [72 \S 3.2] define a system
\[
\colim_{n \to \infty} \; B\cH_\Diff(S^{n-1}) \to \Wh(*) = \Omega^\infty 
\tK(S^0)
\]
of maps to the fiber of the Dennis trace, which becomes a homotopy equivalence
in the limit. It follows that the $K$-theory groups
\[
\colim_{n \to \infty} \; \pi_* B\Diff(\D^n) \; := \; K_*^\cell
\]
(of smooth cell bundles over a point) satisfy
\[
K_i^\cell \otimes \Q \; \cong \; \Q^2 \; \; {\rm if} \; i = 4k > 0
\]
and are zero for other positive $i$, cf eg [42 p 7].\bigskip

\noindent
The resulting parallel manifestations of classical zeta-values in algebraic 
geometry, and in algebraic and differential topology, seem quite remarkable, 
and I am arguing here that they have a unified origin in the fibration
\[
\xymatrix{
\Omega \Wh(*) \ar[r] & B\Diff(\D) \ar[r] & B\Oh }
\]
with odd negative zeta-values originating in the $J$-homomorphism to $Q(S^0)$
on the right, and odd positive zeta-values originating in pseudoisotopy theory
through $K(S^0)$ on the left. The adjoint functors $B$ and $\Omega$ account 
for the shift of homological dimension by two, from $K_{4k-1}(\Z)$ (where 
$\zeta(1-2k)$ lives) to $K_{4k+1}(\Z)$ (where $\zeta(1+2k)$ lives).\bigskip

\noindent
One can hope that this provocative fact might someday provide a 
basis for a theory of {\bf smooth} motives (conceivably involving the
functional equation of the zeta-function), but at the moment even the 
multiplicative structure of $K^\cell_* \otimes \Q$ is obscure to me. 
\bigskip

\noindent
{\bf 3.2.1} Work of Rognes [63], sharpening earlier constructions of
Hatcher [36 \S 6.4], Waldhausen, and B\"okstedt, provides geometrically
motivated generators for $K(S^0)_* \otimes \Q$ by defining a rational 
infinite-loop equivalence
\[
\tw : B(\FO) \to \Wh(*) \;.
\]
Here $\F$ is the monoid of homotopy self-equivalences of the stable sphere 
[56]; I'll write $\fO$ for the spectrum defined by the infinite loopspace 
$\FO$. One of my many debts to this paper's referees is the construction 
of a rational equivalence 
\[
(S^0 \vee \Sigma kO) \otimes \Q \to (S^0 \vee \Sigma \fO) \otimes \Q
\]
of ring-spectra (with simple multiplication) via the zigzag
\[
\xymatrix{
B\Oh = * \times_\Oh E\Oh & \ar[l] F \times_\Oh E\Oh \ar[r] & F \times_\Oh * =
\FO }
\]
of maps of infinite loopspaces; together with Rognes's construction, this 
defines an equivalence 
\[
\w : (S^0 \vee \Sigma kO) \otimes \Q \to K(S^0) \otimes \Q 
\]
of rational ring-spectra (alternately: of DGAs with trivial differentials 
and product structure). \bigskip

\noindent
{\bf Proposition} The resulting homomorphism
\[
\w_* :\Q \oplus kO_*[1] \otimes \Q \to K_*(S^0) \otimes \Q
\]
presents the rationalization of $K(S^0)$ as a square-zero extension of $\Q$
by an ideal
\[
\Q\{\sigma v^k \:|\: k \geq 1 \}
\]
($|v| = 4$) with trivial multiplication.\bigskip

\noindent
{\bf 3.2.2} Writing $S^0[X_+]$ for the suspension spectrum of a space $X$ 
emphasizes the similarity of that construction to the free abelian group
generated by a set. The equivalence
\[
\Maps_{\rm Spaces}(\Omega^\infty Z_+,\Omega^\infty Z_+) = 
\Maps_{\rm Spectra}(S^0[\Omega^\infty Z_+],Z)
\]
sends the identity map on the left to a stabilization morphism 
\[
S^0[\Omega^\infty Z_+] \to Z
\]
of spectra: for example, if $Z = \Sigma kO$ then $\Omega^\infty Z$ is
the Bott space $\SUSO$, and the extension 
\[
S^0[\SUSO_+] \to S^0 \vee \Sigma kO 
\]
of stabilization by the collapse map $\SUSO \to S^0$ to a map of
ring-spectra (with the target regarded as a square-zero extension) 
is the  product-killing quotient 
\[
e_{4k+1} \mapsto \sigma v^k : H_*(\SUSO,\Q) = E(e_{4k+1} \:|\; k \geq 1) 
\otimes \Q \to \Q \oplus \Q\{\sigma v^k \;|\; k \geq 1\} \;. 
\] 
\bigskip

\noindent
{\bf 3.2.3} The K\"unneth spectral sequence 
\[
\Tor^{H_*(\SUSO,\Q)}(H_*(S^0,\Q),H_*(S^0,\Q)) \Rightarrow 
H_*((S^0[\SUSO_+]^\dagger,\Q)
\]
[32 IV \S 4.1] for the rational homology of $S^0 \wedge^L_{S^0[\SUSO_+]} S^0$ 
collapses, yielding an isomorphism of its target with the algebra of symmetric 
functions on generators of degree $4k+2, \; k \geq 0$. It is algebraically
isomorphic to the Rothenberg-Steenrod spectral sequence
\[ 
\Tor^{H_*(\SUSO,\Q)}(\Q,\Q)) \Rightarrow H_*(\SpSU,\Q) \;,
\]
[59 \S 7.4] for $B(\SUSO)$, allowing us to identify $S^0_\Q[\SpSU_+]$ 
with the covariant Koszul dual of $S^0_\Q[\SUSO_+]$. The composition
\[
S^0_\Q[\SpSU_+] = (S^0 \wedge^L_{S^0[\SUSO_+]} S^0)_\Q \to
\]
\[
\to (S^0 \wedge^L_{S^0 \vee \Sigma kO} S^0)_\Q \to (S^0 \wedge_{K(S^0)} 
S^0)_\Q = K(S^0)^\dagger_\Q
\]
represents the abelianization map
\[
\sG(\f^*) \to \sG(\f^*_\ab) \; ( = \Spec \; H_*(\SpSU,\Q))
\]
of \S 2.4 above. \bigskip

\noindent
{\bf 3.2.4 Remarks}\medskip

\noindent
1) $v^2$ is twice the Bott periodicity class. \bigskip

\noindent
2) The arguments above are based on the equivalence, over the rationals,
of $K(S^0)$ and $K(\Z)$. In a way this is analogous to the isomorphism 
between singular (Betti) and algebraic de Rham (Grothendieck) cohomology 
of algebraic varieties. Nori [48 Theorem 6] formulates the theory of periods
in terms of functions on the torsor of isomorphisms between these theories;
from this point of view zeta-values appear as functions on
$\Spec \; (K(\Z)_*K(S^0))$, viewed as a torsor relating arithemetic
geometry to differential topology. \bigskip

\noindent
{\bf 3.3} The Tannakian category of mixed Tate motives over $\Z$ constructed
by Deligne and Goncharov is equivalent to the category of linear 
representations of the motivic group $\GT_\MTM$ of that category (thought
to be closely related to Drinfel'd's prounipotent version of the 
Grothendieck-Teichm\"uller group [2 \S 25.9.4; 28; 73 \S 6.1, Prop 9.1]).
At the end of a later paper Goncharov describes the Hopf algebra 
$\cH_{\GT^*_\MTM)}$ of functions on this motivic group in some detail: in 
particular [35 \S 8.2 Theorem 8.2, \S 8.4 exp (110)] he identifies it as the 
cotensor algebra $T^\bullet(\sK)$, where
\[
\sK := \oplus_{n \geq 1} K_{2n-1}(\Z) \otimes \Q
\]
(regarded as a graded module with $K_{2n-1}$ situated in degree $n$). \bigskip

\noindent
The composition of the pseudo-isotopy map $\w_*$ of \S 3.2.1 with Waldhausen's
isomorphism $K(S^0) \otimes \Q \cong K(\Z) \otimes \Q$ identifies the free
graded Lie algebra on $\sK$ with the free Lie algebra $\f_*$ of \S 2.3iii
above, yielding an isomorphism
\[
\sG(\f_*) \to \GT_\MTM
\]
of proalgebraic groups. Corollary 2.4 then implies the \bigskip

\noindent
{\bf Theorem} The composition
\[
\w_* \circ L\Phi^\dagger : D_b(K_*(S^0) \otimes \Q - \Mod) \to D_b
(\GT_\MTM - \Mod) \;,
\]
defines an equivalence of the homotopy category of rational $K(S^0)$-module 
spectra with the derived category of mixed Tate motives
over $\Z$. \bigskip

\section{Some applications}\bigskip

\noindent
This section discusses some applications of the preceding discussion. 
The first paragraph below is essentially an acknowledgement of ignorance about 
topological cyclic homology. The second discusses some joint work in progress 
[47] with Nitu Kitchloo. The setup and ideas are entirely his; the section 
below sketches how Koszul duality seem to fit in with them. I am indebted 
to Kitchloo for generously sharing these ideas with me. \bigskip

\noindent
The third paragraph summarizes some of the work of Blumberg, Gepner, and 
Tabuada mentioned in the Introduction, concerned with a program for 
constructing enriched decategorifications of their approach to generalized 
motives as small stable $\infty$-categories.
\bigskip

\noindent
{\bf 4.1} Topological cyclic homology [13, 17, \dots] is a powerful tool 
for the study of the algebraic $K$-theory of spaces, and its role in these 
matters deserves discussion here; but at the moment there are technical 
obstructions to telling a coherent story. The current  state of the art 
defines local invariants TC(X;$p$) for a space at each prime $p$ (closely
related to the homotopy quotient of the suspension of the free loopspace
of $X$), whereas the theory of mixed Tate motives over integer rings is 
intrinsically global. For example, the topological cyclic homology of a 
point looks much like the $p$-completion of an ad hoc geometric model
\[
\TC^\geo(S^0) \; \sim \; S^0 \vee \Sigma \CP^\infty_{-1}
\]
[9, 57, 62 \S 3] with
\[
H_*(TC^\geo(S^0),\Z) \cong \Z \oplus \Z\{\sigma t^k \:|\; k \geq -1 \} \;.
\]
The (rational) Koszul dual of this object defines a proalgebraic groupscheme 
associated to a free graded Lie algebra roughly twice as big as $\f^*$, ie 
with generators in topological degree $-2k$ rather than $-2(2k+1)$. A similar 
group appears in work of Connes and Marcolli [25 Prop 5.4] on renormalization 
theory, and topological cyclic homology is plausibly quite relevant to that 
work; but because the global arithmetic properties of topological cyclic 
homology are not yet well understood, it seems premature to speculate further
here; this remark is included only to signal this possible connection to 
physics.\bigskip

\noindent
{\bf 4.2 Example:} Kitchloo [46] has defined a rigid monoidal category 
$\sS$ with symplectic manifolds $(M,\omega)$ as objects, and stable equivalence
classes of oriented Lagrangian correspondences as morphisms. It has a fiber
functor which sends such a manifold (endowed with a compatible almost-complex 
structure) to a Thom spectrum
\[
\sOmega(M) \; = \; \USO(T_M)^{-\zeta}
\]
constructed from the $\USO$-bundle of Lagrangian structures on its stable
tangent space. An Eilenberg-Moore spectral sequence with
\[
E_2 = \Tor^{H^*(BU)}_*(H^*(M),H^*(BSO))
\]
computes $H^*(\USO(T_M))$, and away from the prime two, the equivariant Borel
cohomology 
\[
H^*_\USO(\sOmega(M)) := H^*(\sOmega(M) \times_\USO E(\USO)) 
\]
is naturally isomorphic to $H^*(M)$.\bigskip

\noindent
The functor $\sOmega(-)$ has many of the formal properties of a homology 
theory; for example, when $M$ is a point, $\sOmega := \sOmega(*)$ is a 
ring-spectrum [65], and $\sOmega(-)$ takes values in the category of 
$\sOmega$-modules. Moreover, when $V$ is compact oriented, with the usual 
symplectic structure on its cotangent space,
\[
\sOmega(T^*V) \; \sim \; [V,\sOmega]
\]
[46 \S 2.6] defines a cobordism theory of Lagrangian maps (in the sense of 
Arnol'd) to $V$. \bigskip

\noindent
The composition
\[
\sOmega(M) \to \sOmega(M) \wedge M_+ \to \sOmega(M) \wedge B(\USO)_+
\]
(defined by the map $M \to B(\USO)$ which classifies the bundle 
$\USO(T_M)$ of Lagrangian frames on $M$) makes $\sOmega(M)$ a comodule 
over the Hopf spectrum
\[
\THH(\sOmega) \cong \sOmega \wedge B(\USO)_+ \sim \sOmega[\SpU_+]
\]
(the analog, in this context, of an action of the abelianization $\sG(\f^*_\ab)$ 
[47 \S 4]). The Hopf algebra counit
\[
[1 : S^0_\Q[\USO_+] \to S^0_\Q] \in H^0(\USO,\Q)
\]
provides, via the Thom isomorphism, an augmentation
\[
[\sOmega_\Q \to S^0_\Q] \in H^*(\sOmega,\Q) \cong H^*(\USO,\Q) \;.
\]
{\bf Proposition} The covariant Koszul dual
\[
\sOmega(M)^\dagger_\Q \; := \; \sOmega(M) \wedge^L_{\sOmega_\Q} S^0_\Q
\]
is a comodule over
\[
S^0_\Q \wedge^L_{\sOmega_\Q} S^0_\Q \; \sim \; S^0_\Q \wedge^L_{S^0_\Q[\USO_+]} 
S^0_\Q \; \sim \; S^0[\SpU_+] \;;
\]
by naturality its {\bf contravariant} Koszul dual
\[
R\Hom_{\sOmega_\Q}(S^0_\Q,\sOmega^\dagger_\Q(M)) \cong \sOmega_\Q(M)
\]
inherits an $\sOmega_\Q[\SpU_+]$ - coaction: equivalently, an action of the
abelianized Grothendieck-Teichm\"uller group $\sG(\f^*_\ab)$. $\Box$ \bigskip

\noindent
{\bf Remarks:} \bigskip

\noindent
i) It seems likely that this coaction agrees with the $\THH(\sOmega)$-coaction
described above.\bigskip

\noindent
ii) If $M = T^*V$ is a cotangent bundle, we have an isomorphism
\[
\sOmega^\dagger_\Q(M) \cong H_*(V,\Q) \;.
\]

\noindent
iii) The sketch above is proposed as an analog, in the theory of geometric 
quantization, to work [48 \S 4.6.2, \S 8.4] of Kontsevich on deformation 
quantization. A version of the Grothendieck-Teichm\"uller group acts on the 
Hochschild {\bf co}homology
\[
\HH^*_\C(M) := \Ext^*_{\cO_{M \times M}}(\cO_M,\cO_M) \cong \bigoplus H^*(M,
\Lambda^*T_M))
\]
of a a complex manifold (defined in terms of coherent sheaves of holomorphic
functions on $M \times M$). If $M$ is Calabi-Yau its tangent and cotangent 
bundles can be identified, resulting in an action of the abelianized 
Grothendieck-Teichm\"uller group on the Hodge cohomology of $M$. 
\bigskip

\noindent
Note that $\SpU \sim B\T \times \SpSU$ splits. The action on $\sOmega$ of the 
two-dimensional cohomology class carried by $B\T$ does not seem to come from 
a $K(S^0)^\dagger$ coaction, but rather from variation of the symplectic 
structure. This may be related to Kontsevich's remarks (just after Theorems 7 
and 9) about Euler's constant.
\bigskip

\noindent
{\bf 4.3  Example} Marshalling the forces of higher category theory, Blumberg, 
Gepner, and Tabuada [10] have developed a beautiful approach to the study of 
noncommutative motives, defining symmetric monoidal categories $\mathcal M$ 
(there are several interesting variants [10 \S 6.7, \S 8.10]) whose objects 
are small stable $\infty$-categories (eg of perfect complexes of quasicoherent 
sheaves of modules over a scheme, or of suitably small modules over the 
Spanier-Whitehead dual ring-spectrum $[X,S^0]$ of a finite complex). The 
morphism objects 
\[
\Mor_{\mathcal M}(\cA,\cB)
\]
in these constructions are $K$-theory spectra of categories of exact 
functors between $\cA$ and $\cB$; this defines spectral enrichments over 
the homotopy category of $K(S^0)$-modules [11, Corollary 1.11]. \bigskip

\noindent
The arguments of this paper imply that covariant Koszul duality, as outlined
above, defines versions of these categories with morphism objects 
\[
\Mor_{\HoM}(\cA,\cB)^\dagger_\Q \in K(S^0)^\dagger_\Q - \Comod
\]
which, under suitable finiteness conditions, may be regarded as enriched 
over $D_b(\MTM_\Q(\Z))$. They suggest the existence of categories $\HoM^\dagger$ 
with morphism objects 
\[
\Mor_\HoM(\cA,\cB) \in K(S^0)^\dagger - \Comod
\]
which rationalize to the categories described above. This seems to fit
well with recent work [66] on Konsevich's conjecture on noncommutative
motives over a field [61 \S 4.4]. The theory of cyclotomic spectra [13] 
suggests the existence of related constructions from that point of view,
but (as noted in \S 4.1) their arithmetic properties are not yet very 
well-understood.\bigskip

\noindent
Recently F Brown, using earlier work of Zagier [74], has shown that 
the algebra $\cH_{\GT^*_\MTM}$ is isomorphic to a polynomial algebra 
\[
\Q[\zeta^\mm(w) \:|\: w \in {\rm Lyndon}\{2,3\}] := \Q[\zeta^\mm]
\]
of motivic polyzeta values indexed by certain Lyndon words [cf [20 \S 3 
exp 3.6, \S 8]: working with motivic polyzetas avoids questions of 
algebraic indepence of numerical polyzetas]. This suggests the category 
$\HoM^\dagger_{\Q[\zeta^\mm]}$, with morphism objects
\[
\Mor_\HoM(\cA,\cB)^\dagger_\Q \otimes_{\GT_\MTM} \Q[\zeta^\mm]
\]
as a convenient `untwisted' $\Q[\zeta^\mm]$-linear category of noncommutative 
motives. \bigskip

\bibliographystyle{amsplain}

\end{document}